\documentclass{amsart}

\usepackage{amsmath,amsthm,amsfonts,amscd,amssymb}

\newtheorem{thm}{Theorem}[section]
\newtheorem{cor}[thm]{Corollary}
\newtheorem{lem}[thm]{Lemma}

\theoremstyle{definition}

\theoremstyle{remark}

\begin{document}

\title{Frobenius problem and the covering radius of a lattice}
\author{Lenny Fukshansky \and Sinai Robins}

\address{Department of Mathematics, Mailstop 3368, Texas A\&M University, College Station, Texas 77843-3368}
\email{lenny@math.tamu.edu}
\address{Department of Mathematics, Temple University, Philadelphia, Pennsylvania, 19122}
\email{srobins@math.temple.edu}

\begin{abstract}
Let $N \geq2$ and let $1 < a_1 < \dots < a_N$ be relatively prime integers. Frobenius number of this $N$-tuple is defined to be the largest positive integer that cannot be expressed as $\sum_{i=1}^N a_i x_i$ where $x_1,...,x_N$ are non-negative integers. The condition that $gcd(a_1,...,a_N)=1$ implies that such number exists. The general problem of determining the Frobenius number given $N$ and $a_1,...,a_N$ is NP-hard, but there has been a number of different bounds on the Frobenius number produced by various authors. We use techniques from the geometry of numbers to produce a new bound, relating Frobenius number to the covering radius of the null-lattice of this $N$-tuple. Our bound is particularly interesting in the case when this lattice has equal successive minima, which, as we prove, happens infinitely often.
\end{abstract}

\maketitle

\def\A{{\mathcal A}}
\def\B{{\mathcal B}}
\def\C{{\mathcal C}}
\def\D{{\mathcal D}}
\def\F{{\mathcal F}}
\def\x{{\mathcal H}}
\def\I{{\mathcal I}}
\def\J{{\mathcal J}}
\def\K{{\mathcal K}}
\def\L{{\mathcal L}}
\def\M{{\mathcal M}}
\def\R{{\mathcal R}}
\def\s{{\mathcal S}}
\def\V{{\mathcal V}}
\def\W{{\mathcal W}}
\def\X{{\mathcal X}}
\def\Y{{\mathcal Y}}
\def\H{{\mathcal H}}
\def\aaa{{\mathbb A}}
\def\cee{{\mathbb C}}
\def\Nn{{\mathbb N}}
\def\pee{{\mathbb P}}
\def\que{{\mathbb Q}}
\def\real{{\mathbb R}}
\def\zed{{\mathbb Z}}
\def\gmn{{\mathbb G_m^N}}
\def\qbar{{\overline{\mathbb Q}}}
\def\eps{{\varepsilon}}
\def\vek{{\varepsilon_k}}
\def\ahat{{\hat \alpha}}
\def\bhat{{\hat \beta}}
\def\gt{{\tilde \gamma}}
\def\h{{\tfrac12}}
\def\ba{{\boldsymbol a}}
\def\be{{\boldsymbol e}}
\def\bei{{\boldsymbol e_i}}
\def\bc{{\boldsymbol c}}
\def\bm{{\boldsymbol m}}
\def\bk{{\boldsymbol k}}
\def\bi{{\boldsymbol i}}
\def\bl{{\boldsymbol l}}
\def\bq{{\boldsymbol q}}
\def\bu{{\boldsymbol u}}
\def\bt{{\boldsymbol t}}
\def\bs{{\boldsymbol s}}
\def\bv{{\boldsymbol v}}
\def\bw{{\boldsymbol w}}
\def\bx{{\boldsymbol x}}
\def\bX{{\boldsymbol X}}
\def\bz{{\boldsymbol z}}
\def\bwy{{\boldsymbol y}}
\def\bg{{\boldsymbol g}}
\def\bY{{\boldsymbol Y}}
\def\bL{{\boldsymbol L}}
\def\baa{{\boldsymbol\alpha}}
\def\bb{{\boldsymbol\beta}}
\def\bet{{\boldsymbol\eta}}
\def\bxi{{\boldsymbol\xi}}
\def\bo{{\boldsymbol 0}}
\def\bol{{\boldsymbol 1}_L}
\def\ep{\varepsilon}
\def\p{\boldsymbol\varphi}
\def\q{\boldsymbol\psi}
\def\rank{\operatorname{rank}}
\def\aut{\operatorname{Aut}}
\def\lcm{\operatorname{lcm}}
\def\sgn{\operatorname{sgn}}
\def\spn{\operatorname{span}}
\def\md{\operatorname{mod}}
\def\Norm{\operatorname{Norm}}
\def\dim{\operatorname{dim}}
\def\det{\operatorname{det}}
\def\Vol{\operatorname{Vol}}
\def\rk{\operatorname{rk}}
\def\md{\operatorname{mod}}

\section{Introduction}

Let $N \geq 2$ be an integer and let $a_1,...,a_N$ be positive relatively prime integers. Define the {\it Frobenius number} $\F=\F(a_1,...,a_N)$ of this $N$-tuple to be the largest positive integer that cannot be expressed as $\sum_{i=1}^N a_i x_i$ where $x_1,...,x_N$ are non-negative integers. The condition that $gcd(a_1,...,a_N)=1$ implies that such $\F$ exists. The general problem of determining the Frobenius number given $N$ and $a_1,...,a_N$ is NP-hard. For each fixed $N$, however, it is possible to give a polynomial time algorithm for finding the Frobenius number of a given $N$-tuple (see \cite{kannan}). Since there can be no explicit formula for the Frobenius number, it is interesting to produce upper bounds for it. A large amount of work has been done on this problem. The case of $N=2$ is the only one where an explicit formula, known most likely to Sylvester \cite{sylvester}, is available:
\begin{equation}
\label{sylv}
\F(a_1,a_2) = (a_1-1)(a_2-1)-1.
\end{equation}

In a more general case $N \geq 3$, the bounds on the Frobenius number in the literature are vast. Among many others, they include results by Beck, Diaz, and Robins \cite{beck:diaz:robins} produced with the use of bounds on Fourier-Dedekind sums:
\begin{equation}
\label{bdr}
\F \leq \frac{1}{2} \left( \sqrt{a_1 a_2 a_3 (a_1 + a_2 + a_3)} - a_1 - a_2 - a_3 \right),
\end{equation}
as well as earlier results by Erd\"{o}s and Graham \cite{erdos}
\begin{equation}
\label{erd}
\F \leq 2 a_N \left[ \frac{a_1}{N} \right] - a_1,
\end{equation}
by Selmer \cite{selmer}
\begin{equation}
\label{selm}
\F \leq 2 a_{N-1} \left[ \frac{a_N}{N} \right] - a_N,
\end{equation}
and by Vitek \cite{vitek}
\begin{equation}
\label{vit}
\F \leq \left[ \frac{(a_2-1)(a_N-2)}{2} \right] - 1,
\end{equation}
where $[\ ]$ denotes integer part function. See \cite{beck:diaz:robins} for further bibliography. For comparison, here is a lower bound on $\F$ by Aliev and Gruber \cite{iskander}:
\begin{equation}
\label{isk}
\F > \left((N-1)!\ a_1 \dots a_N \right)^{\frac{1}{N-1}} - \sum_{i=1}^N a_i.
\end{equation}
See \cite{iskander} for more information on lower bounds. The objective of this paper is to produce new upper bounds for the Frobenius number when $N \geq 3$.
\smallskip

In \cite{kannan}, Kannan relates the Frobenius number $\F$ to the covering radius of a certain convex body with respect to a certain lattice. More precisely, let
$$\L = \left\{ \bx \in \zed^{N-1} : \sum_{i=1}^{N-1} a_i x_i \equiv 0\ (\md\ a_N) \right\},$$
and define
$$\s = \left\{ \bx \in \real_{\geq 0}^{N-1} : \sum_{i=1}^{N-1} a_i x_i \leq 1 \right\}.$$
Then Theorem 2.5 of \cite{kannan} states that
\begin{equation}
\label{kannan:id}
\F = \mu(\s,\L) - \sum_{i=1}^N a_i.
\end{equation}
where $\mu(\s,\L)$ is the {\it covering radius} (also known as the {\it inhomogeneous minimum}) of $\s$ with respect to $\L$, namely
\begin{equation}
\label{kcover}
\mu(\s,\L) = \inf \left\{ t \in \real_{>0} : t\s + \L = \real^{N-1} \right\}.
\end{equation}
Identity (\ref{kannan:id}) then suggests that one could produce bounds on $\F$ by bounding $\mu(\s,\L)$. This, however, appears difficult, since the standard techniques for bounding a covering radius only work in the case when the convex body is symmetric with respect to the origin, which is clearly not the case here.
\smallskip

Our approach relates the Frobenius number to a covering radius of a Euclidean ball with respect to a different lattice, which is much easier to estimate. Let $\ba = (a_1,...,a_N) \in \zed_{\geq 0}^N$, with $2 \leq a_1 < a_2 < \dots < a_N$ relatively prime, as above. Let
$$L_{\ba}(\bX) = \sum_{i=1}^N a_i X_i,$$
be the linear form in $N$ variables with coefficients $a_1,...,a_N$, and define the lattice
$$\Lambda_{\ba} = \left\{ \bx \in \zed^N : L_{\ba}(\bx) = 0 \right\}.$$
Let $V_{\ba} = \spn_{\real} \Lambda_{\ba}$, then $V_{\ba}$ is an $(N-1)$-dimensional subspace of $\real^N$ and $\Lambda_{\ba} = V_{\ba} \cap \zed^N$ is a lattice of full rank in $V_{\ba}$. Let $B(R)$ be the $(N-1)$-dimensional closed ball of radius $R>0$ centered at the origin in $V_{\ba}$. Then $\Vol_{N-1}(B(R)) = \omega_{N-1} R^{N-1}$, where
\begin{equation}
\label{unit_vol}
\omega_{N-1} = \Vol_{N-1}(B(1)) = \frac{\pi^{\frac{N-1}{2}}}{\Gamma\left( \frac{N+1}{2} \right)}.
\end{equation}
Define the covering radius of the lattice $\Lambda_{\ba}$ to be
\begin{equation}
\label{cover_def}
R_{\ba} = \inf \left\{ R \in \real_{>0} : B(R) + \Lambda_{\ba} = V_{\ba} \right\}.
\end{equation}
It is not difficult to see that $R_{\ba}$ is the radius of the smallest ball that can be circumscribed around the {\it Voronoi cell} of $\Lambda_{\ba}$, which is defined by
$$\V(\Lambda_{\ba}) = \{ \bwy \in V_{\ba} : \| \bwy \| \leq \| \bwy - \bx \|\ \forall\ \bx \in \Lambda_{\ba} \},$$
where $\|\ \|$ stands for the usual Euclidean norm on vectors. Notice that unlike $\mu(\s,\L)$ of (\ref{kcover}), $R_{\ba}$ is a well understood invariant of the lattice. We will discuss it in further details in section~3. The main result of this paper is the following theorem.

\begin{thm} \label{main} Let $N \geq 3$ and let $2 \leq a_1 < a_2 < \dots < a_N$ be relatively prime integers. Write $\ba = (a_1,...,a_N)$, and let $\F = \F(\ba)$ be the Frobenius number of this $N$-tuple. Then
\begin{equation}
\label{bound}
\F \leq \left[ \frac{(N-1) R_{\ba}}{\|\ba\|} \sum_{i=1}^N a_i \sqrt{ \|\ba\|^2 - a_i^2} + 1 \right],
\end{equation}
where $R_{\ba}$ is as in (\ref{cover_def}).
\end{thm}

Our approach uses some classical results from the geometry of numbers. Here is a brief outline of our argument. Let $t$ be a positive integer, and consider the hyperplane in $\real^N$ defined by the equation 
\begin{equation}
\label{def}
\sum_{i=1}^N a_i X_i = t. 
\end{equation}
The intersection of this hyperplane with the positive orthant $\real^N_{\geq 0}$ is an $(N-1)$-dimensonal simplex, call it $S(t)$. An integral point in this simplex corresponds to a solution of (\ref{def}) in non-negative integers, hence for every $t > \F$ such a point must always exist. Moreover, $\F$ is precisely the smallest positive integer such that for each integer $t > \F$ the simplex $S(t)$ contains a point of $\zed^N$. By definition of $R_{\ba}$, a ball of radius $\geq R_{\ba}$ must contain an integer lattice point. On the other hand, it is possible to bound the inradius of the simplex $S(t)$ from below using a standard isoperimetric inequality. Combining these two estimates produces a value $t_*$ large enough so that for every $t \geq t_*$ the simplex $S(t)$ is guaranteed to contain an integral point.
\smallskip

A particularly nice explicit bound for $\F$ can be derived from Theorem \ref{main} for a special class of latices $\Lambda_{\ba}$. For each $1 \leq i \leq N-1$, the $i$-th {\it successive minimum} $\lambda_i$ of $\Lambda_{\ba}$ is defined to be the infimum of all $\lambda > 0$ such that $B(\lambda) \cap \Lambda_{\ba}$ contains $i$ non-zero linearly independent vectors in $V_{\ba}$. Hence $1 \leq \lambda_1 \leq ... \leq \lambda_{N-1}$. If $\lambda_1 = \dots = \lambda_{N-1}$, we say that $\Lambda_{\ba}$ is an {\it ESM lattice} (equal successive minima). This is a very important class of lattices, which are widely used for instance in coding theory (see \cite{esm}).

\begin{cor} \label{esm_case} Let the notation be as above. Then
\begin{equation}
\label{esm_bound}
\F \leq \left[ \frac{\lambda_{N-1} (N-1)^2 \sum_{i=1}^N a_i \sqrt{ \|\ba\|^2 - a_i^2}}{\lambda_1 (\|\ba\|^{N-2} \omega_{N-1})^{\frac{1}{N-1}}} + 1 \right],
\end{equation}
where $\omega_{N-1}$ is as in (\ref{unit_vol}). In case $\Lambda_{\ba}$ is an ESM lattice, $\lambda_{N-1} = \lambda_1$ in (\ref{esm_bound}).
\end{cor}

\noindent
One interesting feature of our bounds (\ref{bound}) and (\ref{esm_bound}) is that they depend symmetrically on all numbers $a_1, \dots, a_N$, unlike the previously known bounds (\ref{bdr})~-~(\ref{vit}).
\smallskip

In section~2 of this paper we prove Theorem \ref{main}. In section~3 we discuss the ESM case, deriving Corollary \ref{esm_case}, as well as some other related cases using additional tools from the classical geometry of numbers. We also show some examples and exhibit some computational data comparing our bounds to the previously known ones quoted in (\ref{bdr}) - (\ref{vit}). In particular, when $\Lambda_{\ba}$ is an ESM lattice, Corollary \ref{esm_case} will often produce a better bound on $\F$ than (\ref{bdr}) - (\ref{vit}). We discuss this further in section~3. In section~4 we prove that $\Lambda_{\ba}$ is an ESM lattice for infinitely many $N$-tuples $\ba$. In fact, in Theorem \ref{esm-4} we construct an explicit infinite family of ESM lattices $\Lambda_{\ba}$ parametrized by integer values of a single variable $t$ when $N=4$. We also explain how families like this can be constructed in higher dimensions. Finally we demonstrate that for all such infinite families of ESM lattices $\Lambda_{\ba}$ our bound (\ref{esm_bound}) on $\F(\ba)$ is significantly better than the previously known ones.
\smallskip

\section{Proof of Theorem \ref{main}}

Let the notation be as in section~1 above. For each $t \in \zed_{\geq 0}$ consider the hyperplane lattice
$$\Lambda_{\ba}(t) = \left\{ \bx \in \zed^N : L_{\ba}(\bx) = t \right\},$$
and let $V_{\ba}(t) = \spn_{\real} \Lambda_{\ba}(t)$ be the corresponding hyperplane. Fix $\bu_t \in \Lambda_{\ba}(t)$, and define a translation $f_t : V_{\ba} \rightarrow V_{\ba}(t)$ given by $f_t(\bx) = \bx + \bu_t$ for each $\bx \in V_{\ba}$. Then $f_t$ is bijective and preserves distance; moreover, it maps $\Lambda_{\ba}$ bijectively onto $\Lambda_{\ba}(t)$. 

Notice that $S(t) = V_{\ba}(t) \cap \real_{\geq 0}^N$ is an $(N-1)$-dimensional simplex in $\real^N$ with vertices $\bv_i = \frac{t}{a_i} \be_i$ for each $1 \leq i \leq N$, where $\be_1,...,\be_N$ are the standard basis vectors. For each $2 \leq i \leq N$ define
$$\bw_i = (\bv_i - \bv_1)^T = \left( - \frac{t}{a_1}, 0, ..., 0, \frac{t}{a_i}, 0, ..., 0 \right),$$
and let $W$ be the $(N-1) \times N$ matrix with row vectors $\bw_2,...,\bw_N$. By Gram determinant formula
\begin{equation}
\label{sv1}
\Vol_{N-1}(S(t)) = \frac{\sqrt{\det(W W^T)}}{(N-1)!} 
\end{equation}
It is easy to see that 
$$WW^T = \frac{t^2}{a_1^2} \left( \begin{matrix}
\frac{a_1^2+a_2^2}{a_2^2}&1&\hdots&1\\
1&\frac{a_1^2+a_3^2}{a_3^2}&\hdots&1\\
\vdots&\vdots&\ddots&\vdots\\
1&1&\hdots&\frac{a_1^2+a_N^2}{a_N^2}
\end{matrix} \right),$$
is an $(N-1) \times (N-1)$ symmetric matrix. We want to compute $\det(WW^T)$. For this we will need the following lemma.

\begin{lem} \label{det_comp} Let 
$$\A = \left( \begin{matrix}
\alpha_1&1&\hdots&1\\
1&\alpha_2&\hdots&1\\
\vdots&\vdots&\ddots&\vdots\\
1&1&\hdots&\alpha_k
\end{matrix} \right),$$
be a $k \times k$ symmetric matrix, $k \geq 2$. Then
\begin{equation}
\label{det_formula}
\det(\A) = \prod_{i=1}^k (\alpha_i-1) + \sum_{i=1}^k \left\{ \prod_{j=1,\ j \neq i}^k (\alpha_j-1) \right\}.
\end{equation}
\end{lem}

\proof
It is easy to notice that $\det(\A) = \det(\B)$, where
$$\B = \det \left( \begin{matrix}
\alpha_1-1&0&\hdots&0&1-\alpha_k\\
0&\alpha_2-1&\hdots&0&1-\alpha_k\\
\vdots&\vdots&\ddots&\vdots&\vdots\\
0&0&\hdots&\alpha_{k-1}-1&1-\alpha_k\\
1&1&\hdots&1&\alpha_k
\end{matrix} \right).$$
We will prove identity (\ref{det_formula}) for $\det(\B)$ by induction on $k$. If $k=2$, then
$$\det(\B) = \det \left( \begin{matrix}
\alpha_1-1&1-\alpha_2\\
1&\alpha_2\\
\end{matrix} \right) = (\alpha_1-1) \alpha_2 + (\alpha_2-1),$$
which is (\ref{det_formula}). Assume $k>2$. Then, by Laplace's expansion combined with the induction hypothesis, we obtain
\begin{eqnarray*}
\det(\B) & = & (\alpha_1-1) \det \left( \begin{matrix}
\alpha_2-1&0&\hdots&0&1-\alpha_k\\
0&\alpha_3-1&\hdots&0&1-\alpha_k\\
\vdots&\vdots&\ddots&\vdots&\vdots\\
0&0&\hdots&\alpha_{k-1}-1&1-\alpha_k\\
1&1&\hdots&1&\alpha_k
\end{matrix} \right) \\
& + & (-1)^{k+1} \det \left( \begin{matrix}
0&0&\hdots&0&1-\alpha_k\\
\alpha_2-1&0&\hdots&0&1-\alpha_k\\
\vdots&\vdots&\ddots&\vdots&\vdots\\
0&0&\hdots&0&1-\alpha_k\\
0&0&\hdots&\alpha_{k-1}-1&1-\alpha_k\\
\end{matrix} \right) \\
& = & (\alpha_1-1) \left( \prod_{i=2}^k (\alpha_i-1) + \sum_{i=2}^k \left\{ \prod_{j=2,\ j \neq i}^k (\alpha_j-1) \right\} \right) \\
& + & (-1)^{k+1+k} (1-\alpha_k) \det \left( \begin{matrix}
\alpha_2-1&\hdots&0\\
\vdots&\ddots&\vdots\\
0&\hdots&\alpha_{k-1}-1\\
\end{matrix} \right) \\
& = & \prod_{i=1}^k (\alpha_i-1) + \sum_{i=2}^k \left\{ \prod_{j=1,\ j \neq i}^k (\alpha_j-1) \right\} + \prod_{j=2}^k (\alpha_j-1).
\end{eqnarray*}
This completes the proof.
\endproof

Applying Lemma \ref{det_comp} to $WW^T$, a direct computation shows that
\begin{equation}
\label{sv2}
\det(WW^T) = \frac{t^{2(N-1)} \|\ba\|^2}{\prod_{i=1}^N a_i^2},
\end{equation}
and so combining (\ref{sv1}) with (\ref{sv2}) we obtain
\begin{equation}
\label{simp_vol}
\Vol_{N-1}(S(t)) = \frac{t^{N-1} \|\ba\|}{(N-1)! \prod_{i=1}^N a_i}.
\end{equation}
We also need to compute the surface area $A_{N-1}(S(t))$. Notice that $S(t)$ has $N$ faces $F_1(t),...,F_N(t)$ with each $F_i(t)$ being an $(N-2)$-dimensional simplex with vertices $\bv_1,...,\bv_{i-1},\bv_{i+1},...,\bv_N$. Then, applying (\ref{simp_vol}) in one less dimension we see that for each $1 \leq i \leq N$. 
$$\Vol_{N-2}(F_i(t)) = \frac{t^{N-2} \|\baa_i\|}{(N-2)! \prod_{j=1,\ j \neq i}^N a_j},$$
where $\baa_i = (a_1,...,a_{i-1},a_{i+1},...,a_N)$. Then
\begin{equation}
\label{simp_area}
A_{N-1}(S(t)) = \sum_{i=1}^N \Vol_{N-2}(F_i(t)) = \frac{t^{N-2} \sum_{i=1}^N \|\baa_i\| a_i}{(N-2)! \prod_{i=1}^N a_i}.
\end{equation}
Write $r(t)$ for the inradius of $S(t)$, i.e. the radius of the largest ball that can be inscribed into $S(t)$. By a standard isoperimetric inequality for the inradius of a simplex (see for instance (9) of \cite{inradius})
\begin{equation}
\label{inr}
r(t) \geq \frac{\Vol_{N-1}(S(t))}{A_{N-1}(S(t))} = \frac{t \|\ba\|}{(N-1) \sum_{i=1}^N \|\baa_i\| a_i},
\end{equation}
where the last identity follows by combining (\ref{simp_vol}) and (\ref{simp_area}). Let us choose a positive integer $t$ such that $r(t) \geq R_{\ba}$. By (\ref{inr}) we see that it suffices to take 
\begin{equation}
\label{t_bound}
t = \left[ \frac{(N-1) R_{\ba}}{\|\ba\|} \sum_{i=1}^N \|\baa_i\| a_i + 1 \right].
\end{equation}
We will write $t_*$ for the value of $t$ as in (\ref{t_bound}). Let $t \geq t_*$, and let $B_t(r(t))$ be the $(N-1)$-dimensional ball of radius $r(t)$ contained in $S(t)$. Then $f_t^{-1}(B_t(r(t)))$ is an $(N-1)$-dimensional ball of radius $r(t) \geq R_{\ba}$ in $V_{\ba}$. By definition of $R_{\ba}$ in (\ref{cover_def}), we see that whenever $R \geq R_{\ba}$ the translated ball $B(R)+\bx$ will contain at least one nonzero lattice point for every $\bx \in V_{\ba}$, and hence $f_t^{-1}(B_t(r(t)))$ contains a nonzero point of $\Lambda_{\ba}$. Therefore $B_t(r(t))$ contains a point of $\Lambda_{\ba}(t)$, that is $\Lambda_{\ba}(t) \cap \zed^N_{\geq 0}$ is not empty for each integer $t \geq t_*$. Therefore $\F \leq t_*$, and observing that $\|\baa_i\| = \sqrt{ \|\ba\|^2 - a_i^2}$ for each $1 \leq i \leq N$ finishes the proof.
\smallskip

{\it Remark.} It is possible to replace (\ref{inr}) by stronger versions of this isoperimetric inequality, which follow from the proof of Wills conjecture and its various strengthenings (see, for instance, (4), (6), and Theorem 4 of \cite{wills}). This may lead to a slightly better although much less readable bound than (\ref{bound}). 
\smallskip

\section{Corollaries}

In this section we discuss consequences of Theorem \ref{main}, in particular we derive Corollary \ref{esm_case}. Let $N \geq3$, and let all the notation be as in sections~1 and~2 above. First of all notice that if for some $1 < i \leq N$ we can express $a_i$ in the form
\begin{equation}
\label{bad}
a_i = \sum_{j=1}^{i-1} a_j x_j,
\end{equation}
for some nonnegative integers $x_1,...,x_{i-1}$, then 
\begin{equation}
\label{non_reduced}
\F(a_1,...,a_N)=\F(a_1,...,a_{i-1},a_{i+1},...,a_N).
\end{equation}
We will call the relatively prime $N$-tuple $\ba$ {\it reduced} if (\ref{bad}) is not true for any $i$. By (\ref{non_reduced}), every relatively prime $N$-tuple can be reduced to a relatively prime reduced $k$-tuple for some $1 \leq k \leq N$ by eliminating all $a_i$'s for which (\ref{bad}) is true. Moreover, if $a_1=2$, then there must exist $1 < i \leq N$ such that $a_i$ is odd, since $gcd(a_1,...,a_N)=1$; let $i$ be the smallest such index. It is easy to see that in this case $\F = a_i-1$. In particular, if $\ba$ is reduced, then $\F = a_2-1$. Hence we can conclude that either $a_1 \geq 3$, or
\begin{equation}
\label{not3}
\F \leq a_N-1.
\end{equation}
From here on we will assume that $\ba$ is reduced and $a_1 \geq3$.
\smallskip

Fix a basis $\bx_1,...,\bx_{N-1}$ for $\Lambda_{\ba}$ in $\real^N$, and write $X = (\bx_1\ ...\ \bx_{N-1})$ for the corresponding $N \times (N-1)$ basis matrix. Let $\I$ be the collection of all subsets $I$ of $\{1,...,N\}$ of cardinality $(N-1)$. For each $I \in \I$ let $I'$ be its complement, i.e. $I' = \{1,...,N\} \setminus I$. Clearly $|\I| = N-1$. For each $I \in \I$, write $X_I$ for the $(N-1) \times (N-1)$ submatrix of $X$ consisting of all those rows of $X$ which are indexed by $I$, and $a_{I'}$ for the coordinate of $\ba$ indexed by $I'$. By the duality principle of Brill-Gordan \cite{gordan:1} (also see Theorem 1 on p. 294 of \cite{hodge:pedoe})
\begin{equation}
\label{duality}
\det (X_I) = (-1)^{N+1-I'} a_{I'}.
\end{equation}
Therefore coordinates of $\ba$ can be thought of as {\it Grassmann coordinates} of $\Lambda_{\ba}$ up to $\pm$ signs (some sources also call them {\it Plucker coordinates}). They are well defined in the sense that they do not depend on the choice of the basis (see \cite{hodge:pedoe} for details). Then, by the Cauchy-Binet formula (see for instance \cite{matrices})
\begin{equation}
\label{det}
\det(\Lambda_{\ba}) = \sqrt{\det(X X^t)} = \| \ba \|.
\end{equation}

Let $\lambda_1,...,\lambda_{N-1}$ be the successive minima for $\Lambda_{\ba}$ as defined in section~1. An immediate observation is that since $\ba$ is reduced,
\begin{equation}
\label{slbound}
2 \leq \lambda_1 \leq \dots \leq \lambda_{N-1}.
\end{equation}
Indeed, if $\lambda_1 < 2$, then there must exist $\boldsymbol 0 \neq \bx \in \Lambda_{\ba}$ with $\|\bx\| < 2$, hence at most three of its coordinates are non-zero, call them $x_i,x_j,x_k$, $1 \leq i < j < k \leq N$. Assume $x_i \geq 0$ (take $-\bx$ otherwise). Then either $x_i,x_j=1$ and $x_k=-1$, or one of them is $0$ and the other two are $\pm 1$ and $\mp 2$ respectively. In the first case it must therefore be that $a_k=a_i+a_j$ while the second case implies that one of the coordinates of $\ba$ is a multiple of another. Both of these conclusions contradict the assumption that $\ba$ is reduced.

Combining Minkowski's second convex body theorem (see \cite{cass:geom}, p. 203) with (\ref{det}), we obtain
\begin{equation}
\label{mink}
\lambda_1 \dots \lambda_{N-1} \leq \frac{2^{N-1} \|\ba\|}{\omega_{N-1}}.
\end{equation}
Combining Jarnik's inequality (see Theorem 1 on p. 99 of \cite{lek}) with (\ref{mink}), we obtain a bound on $R_{\ba}$:
\begin{equation}
\label{inh}
R_{\ba} \leq \frac{1}{2} \sum_{i=1}^{N-1} \lambda_i \leq \frac{N-1}{2} \lambda_{N-1} \leq \frac{2^{N-2} (N-1) \|\ba\|}{\omega_{N-1} \lambda_1 \dots \lambda_{N-2}}.
\end{equation}
Then Theorem \ref{main} combined with (\ref{slbound}) and (\ref{inh}) yields a general bound
\begin{equation}
\label{gb}
\F \leq \left[ \frac{(N-1)^2}{\omega_{N-1}} \sum_{i=1}^N a_i \sqrt{ \|\ba\|^2 - a_i^2} + 1 \right],
\end{equation}
however we can do much better for more specialized classes of lattices $\Lambda_{\ba}$. Combining (\ref{mink}) and (\ref{inh}), we obtain
\begin{equation}
\label{ratio}
R_{\ba} \leq \frac{\lambda_1}{2} \sum_{i=1}^{N-1} \frac{\lambda_i}{\lambda_1} \leq \lambda_1 \frac{(N-1)\lambda_{N-1}}{2 \lambda_1} \leq \frac{(N-1) \lambda_{N-1}}{\lambda_1} \left( \frac{\|\ba\|}{\omega_{N-1}} \right)^{\frac{1}{N-1}},
\end{equation}
which, combined with Theorem \ref{main}, immediately implies Corollary \ref{esm_case}. Clearly the bound of Corollary \ref{esm_case} becomes better when the ratio $\frac{\lambda_{N-1}}{\lambda_1}$ is small, and especially in case $\Lambda_{\ba}$ is an ESM lattice.
\bigskip

We will now show a few examples of $\ba$ such that $\Lambda_{\ba}$ is an ESM lattice for which (\ref{esm_bound}) of Corollary \ref{esm_case} produces a better bound on the Frobenius number than (\ref{bdr}) - (\ref{vit}). In the following comparison tables of the bounds (\ref{bdr}) - (\ref{vit}) with (\ref{esm_bound}), $\lambda_{\ba}$ stands for the common value of the successive minima of $\Lambda_{\ba}$. First let $N=4$.
\begin{center} 
\begin{tabular}{|l|l|l|l|} \hline
{\em 4-tuple $\ba$} & {\em $\lambda_{\ba}$} & {\em min (\ref{bdr}) - (\ref{vit})} & {\em (\ref{esm_bound})} \\ \hline \hline
9337, 9961, 11593, 67367 & $\sqrt{1802}$ & 91235853 (\ref{bdr}) & 10995433 \\ \hline
33199, 38351, 47759, 152057 & $\sqrt{3218}$ & 1346684400 (\ref{bdr}) & 55055950 \\ \hline
\end{tabular}
\end{center}
Next let $N=5$.
\begin{center} 
\begin{tabular}{|l|l|l|l|} \hline
{\em 5-tuple $\ba$} & {\em $\lambda_{\ba}$} & {\em min (\ref{bdr}) - (\ref{vit})} & {\em (\ref{esm_bound})} \\ \hline \hline
39221, 46967, 47869, & & & \\
\ \ \ \ \ \ \ \ \ 62839, 206749 & $\sqrt{524}$ & 1719019240 (\ref{bdr}) & 66231577 \\ \hline
1867558, 2348176, 2918749, & & & \\ 
\ \ \ \ \ \ \ \ \ \ \ \ 5249843, 26695349 & $\sqrt{5591}$ & 4778060891200 (\ref{bdr}) & 14595157176 \\ \hline
\end{tabular}
\end{center}
Finally let $N=6$.
\begin{center} 
\begin{tabular}{|l|l|l|l|} \hline
{\em 6-tuple $\ba$} & {\em $\lambda_{\ba}$} & {\em min (\ref{bdr}) - (\ref{vit})} & {\em (\ref{esm_bound})} \\ \hline \hline
6595, 90709, 110483, & & & \\ 
\ 121833, 147472, 462217 & $\sqrt{209}$ & 1015946371 (\ref{erd}) & 168600688 \\ \hline
5958323, 14864655, & & & \\
\ 19945128, 28191201, & & & \\
\ 28507523, 117697394 & $\sqrt{1915}$ & 134180083643479 (\ref{bdr}) & 104669816535 \\ \hline
\end{tabular}
\end{center}

It is of course possible to come up with numerous such examples for these and higher dimensions. In fact, in the next section we will show that $\Lambda_{\ba}$ is an ESM lattice for infinitely many $\ba$.

\section{ESM Lattices}

Let $N \geq 4$. In this section we will describe a procedure that allows to construct infinite families of sublattices of $\zed^N$ of rank $N-1$ which have equal successive minima and are of the form $\Lambda_{\ba}$ for $N$-tuples $\ba$ of relatively prime positive integers $1 < a_1 < \dots < a_N$. 

We start with some additional notation, following \cite{near:ort}. An ordered collection of linearly independent vectors $\{ \bx_1, \dots, \bx_k \} \subset \zed^N$, $2 \leq k \leq N$, is called {\it nearly orthogonal} if for each $1 < i \leq k$ the angle between $\bx_i$ and the subspace of $\real^N$ spanned by $\bx_1, \dots, \bx_{i-1}$ is in the interval $\left[ \frac{\pi}{3}, \frac{\pi}{2} \right]$. In other words, this condition means that for each $1 < i \leq k$
\begin{equation}
\label{near}
\frac{| < \bx_i, \bwy > |}{\|\bx_i\| \|\bwy\|} \leq \frac{1}{2},
\end{equation}
for all non-zero vectors $\bwy \in \spn_{\real} \{ \bx_1, \dots, \bx_{i-1} \}$, where $<\ ,\ >$ stands for the usual inner product on $\real^N$. The following result is Theorem 1 of \cite{near:ort}; it is our main tool in this section.

\begin{thm} [\cite{near:ort}] \label{no} Suppose that an ordered basis $\{ \bx_1, \dots, \bx_k \}$ for sublattice $\Lambda$ of $\zed^N$ of rank $1 < k \leq N$ is nearly orthogonal. Then it contains the shortest non-zero vector of $\Lambda$.
\end{thm}

\noindent
In particular, if all vectors $\bx_1, \dots, \bx_k$ of Theorem \ref{no} have the same norm, then $\Lambda$ is an ESM lattice. We are now ready to describe our construction for infinite families of ESM lattices.
\smallskip

Let $\bx_1 = (t_1, \dots ,t_N)$ be a variable vector, and write $S_N$ for the symmetric group on $N$ letters where $id$ stands for the identity permutaion. Assume that there exist $id = \sigma_1, \sigma_2, \dots, \sigma_{N-1} \in S_N$ and $N(N-1)$ integers $m_{11}, \dots, m_{(N-1)N} \in \{0,1\}$ such that 
$$\bx_i = \left( (-1)^{m_{i1}} t_{\sigma_i(1)}, \dots, (-1)^{m_{iN}} t_{\sigma_i(N)} \right),\ 1 \leq i \leq N-1,$$
satisfy the following conditions for infinitely many positive integer values of the variables $t_1, \dots, t_N$:

\begin{enumerate}
\item $\bx_1, \dots, \bx_{N-1}$ are linearly independent,
\item For each $1 \leq i \leq N$ the corresponding Grassmann coordinate $\det(X_{I_i})$ of the matrix $X = (\bx_1 \dots \bx_{N-1})^t$ satisfies the condition 
$$(-1)^{N+1-i} \det(X_{I_i}) > 0,$$
where $I_i=\{1, \dots, N\} \setminus \{i\}$,
\item Absolute values of Grassmann coordinates of $X$ are relatively prime,
\item $\{ \bx_1, \dots, \bx_{N-1} \}$ is a nearly orthogonal collection of vectors.
\end{enumerate}
\smallskip

\noindent
Then, by construction and by Theorem \ref{no}, for each such $N$-tuple $t_1, \dots, t_N$ the lattice 
$$\spn_{\zed} \{\bx_1, \dots, \bx_{N-1}\}$$
is ESM and of the form $\Lambda_{\ba}$ where $\ba$ is the vector with coordinates 
$$a_i = (-1)^{N+1-i} \det(X_{I_i}),$$ 
for each $1 \leq i \leq N$; the last statement follows by (\ref{duality}). This would mean that there exist infinite families of ESM lattices of the form $\Lambda_{\ba}$. It appears to be possible to find such permutations for each $N$. As an example, we exhibit such a family for $N=4$.

\begin{thm} \label{esm-4} Let $t \in \zed_{> 0}$, and define
\begin{eqnarray}
\label{at}
a_1(t) = 6t^2 -13t - 216,\ a_2(t) = 6t^2 - 125,\nonumber \\
a_3(t) = 7t^2 - 174,\ a_4(t) = t^3 - 36t - 78.
\end{eqnarray}
Then for each $t \in \zed_{> 0}$, $\ba(t) = (a_1(t),a_2(t),a_3(t),a_4(t)) \in \zed^4$, and there exist infinitely many positive integer values of $t$ such that
\begin{equation}
\label{pos} 
0 < a_1(t) < a_2(t) < a_3(t) < a_4(t),
\end{equation}
\begin{equation}
\label{gcd}
gcd \left(a_1(t),a_2(t),a_3(t),a_4(t) \right) = 1,
\end{equation}
and the lattice 
$$\Lambda_{\ba(t)} = \left\{ \bx \in \zed^4 : \sum_{i=1}^4 a_i(t) x_i = 0 \right\}$$
is ESM. Moreover, for each such $\ba(t)$ the minimum of bounds (\ref{bdr})~-~(\ref{vit}) on the Frobenius number $\F(\ba(t))$ is $O(t^4)$ while our bound (\ref{esm_bound}) is $O(t^3)$. For instance, $\ba(t)$ has these properties for all $t = 13s+2$, where $s \geq 2$ is an integer.
\end{thm}

\proof
Let $t \in \zed_{> 0}$ and define
\begin{equation}
\label{circ_bas}
\bx_1(t) = (-7, t, 6, -6),\ \bx_2(t) = (-6, 7, t, -6),\ \bx_3(t) = (-6, -6, 7, t).
\end{equation}
A direct computation shows that
$$\Lambda_{\ba(t)} = \left\{ \bx \in \zed^4 : \sum_{i=1}^4 a_i(t) x_i = 0 \right\} = \spn_{\zed} \{\bx_1(t),\bx_2(t),\bx_3(t)\},$$
where $\ba(t)$ is as in (\ref{at}), and its coordinates can be seen to have no common roots. In particular, $\Lambda_{\ba(t)}$ has rank $3$ and basis vectors $\bx_1(t),\bx_2(t),\bx_3(t)$ are linearly independent for all real values of $t$. Also notice that for each $t \geq 10$, (\ref{pos}) is satisfied. 

To demonstrate that (\ref{gcd}) holds infinitely often, notice that
$$gcd(a_1(t),a_2(t),a_3(t),a_4(t)) \leq gcd(a_2(t),a_3(t)),$$
and define $d(t) = gcd(a_2(t),a_3(t)) = gcd(7t^2-174, 6t^2-125)$. Then $d(t)$ must divide both
$$a_2(t) - a_3(t) = t^2-49, \ 7a_2(t) - 6a_3(t) = 13^2.$$
Notice that if, for instance, $t = 13s+2$ for any $s \in \zed_{>0}$, then
$$t^2-49 = 169s^2+52s-45 \equiv 7\ (\md 13),$$
hence $gcd(t^2-49, 13^2)=1$, and so $d(t)=1$ for all such $t$. This proves that (\ref{gcd}) holds for infinitely many $t \in \zed_{>0}$.

We now want to show that $\{\bx_1(t),\bx_2(t),\bx_3(t)\}$ is a nearly orthogonal ordered collection of vectors for infinitely many $t \in \zed_{>0}$. For this we refer to criterion (\ref{near}) and first observe that
$$\frac{| < \bx_1, \bx_2 > |}{\|\bx_1\| \|\bx_2\|} = \frac{13t + 78}{t^2 + 121} \leq \frac{1}{2},$$
for all $t \geq 28$. Also, for each non-zero vector $\bwy = u\bx_1 + v\bx_2 \in \spn_{\real} \{\bx_1,\bx_2\}$ define
\begin{equation}
\label{cos}
f(u,v) = \frac{t(12u-v)-6(14u-v)}{\sqrt{(t^2+121)\left\{ (u^2+v^2)(t^2+121) + 26uv(t+6) \right\}}} = \frac{< \bx_3, \bwy >}{\|\bx_3\| \|\bwy\|}.
\end{equation}
A computation of the critical points of $f(u,v)$ in Maple shows that if $t \geq 17$ then $-\frac{1}{2} \leq f(u,v) \leq \frac{1}{2}$ for all $u,v \in \real$, not both zero. Hence by criterion (\ref{near}) we conclude that $\{\bx_1(t),\bx_2(t),\bx_3(t)\}$ is a nearly orthogonal ordered collection of vectors for all integers $t \geq 28$. Therefore, by Theorem \ref{no} and remark after it the lattice $\Lambda_{\ba(t)}$ is ESM for all such values of $t$. 

Finally, a direct computation shows that for each $\ba(t)$ as in (\ref{at}) the minimum of bounds (\ref{bdr}) - (\ref{vit}) on the Frobenius number $\F(\ba(t))$ is $O(t^4)$ while bound (\ref{esm_bound}) is $O(t^3)$.

Combining all these observations, we conclude that the statement of the theorem is true for instance for all $t$ of the form 
\begin{equation}
\label{t}
t = 13s+2,
\end{equation}
where $s \geq 2$ is an integer. This completes the proof.
\endproof

Notice in particular that the first example from the table in case $N=4$ in section~3 is precisely of the form (\ref{at}) where $t$ is as in (\ref{t}) with $s=3$. A good strategy to obtain one-parameter infinite families of ESM lattices of the form $\Lambda_{\ba}$ in different dimensions seems to be by a variation on a circulant basis matrix with $\pm$ signs as in (\ref{circ_bas}). In fact, the rest of the examples in the table of section~3 can also be seen to come from such infinite families. 

Moreover, one can see that for a general $N$ if a lattice $\Lambda_{\ba(t)}$ is ESM and is generated by an $(N-1) \times N$ circulant basis matrix with $\pm$ signs similar to (\ref{circ_bas}), call this matrix $X(t)$, then $t$ appears precisely once in every row of $X(t)$ and in all, except for one, columns of $X(t)$. This means that all, except for one, Grassmann coordinates of $X(t)$ in general will be polynomials of degree $N-2$ in $t$, and one will be a polynomial of degree $N-1$. It is not difficult to see that in general in this case the minimum of bounds (\ref{bdr}) - (\ref{vit}) on the Frobenius number $\F(\ba(t))$ will be $O\left(t^{2(N-2)}\right)$ while our bound (\ref{esm_bound}) will be $O\left(t^{N-1} \right)$.

\bibliographystyle{plain}  
\bibliography{frobenius}   

\end{document}